\documentclass{article}
\usepackage{amsmath,amssymb,amsthm}
\usepackage{graphics}

\newtheorem{theorem}{Theorem}[section]
\newtheorem{lemma}{Lemma}[section]
\newtheorem{remark}{Remark}[section]

\newcommand{\tr}{\mathop{\rm Tr}\nolimits}
\newcommand{\rk}{\mathop{\rm rank}\nolimits}
\newcommand{\card}{\mathop{\rm card}\nolimits}
\newcommand{\psd}{\mathop{\rm PD}\nolimits}

\newtheorem{cor}{Corollary}[section]
\newtheorem{prop}{Proposition}[section]
\newtheorem{defn}{Definition}[section]

\title {Multivariate positive definite functions on spheres}

\author {Oleg R. Musin \thanks{University of Texas at Brownsville. Email: oleg.musin@utb.edu}}

\begin{document}
\date{}
\maketitle

\begin{abstract} In 1942  I. J. Schoenberg  proved that a function  is positive definite  in the unit sphere if and only if this function is a positive  linear combination of the Gegenbauer  polynomials.  In this paper we extend Schoenberg's theorem  for multivariate
 Gegenbauer polynomials. This extension derives new positive semidefinite constraints for the distance distribution which can be applied for spherical codes.
\end {abstract}

\section{Introduction}

Let $M$ be a metric space with a distance function $\tau.$ A real continuous function $f(t)$  
is said to be positive definite (p.d.) in $M$ if
for arbitrary points $p_1,\ldots,p_r$ in $M$, real variables $x_1,\ldots,x_r$, and arbitrary $r$  we have
$$
\sum\limits_{i=1}^r\sum\limits_{j=1}^r {f(t_{ij})\,x_ix_j}\ge 0, \quad t_{ij}=\tau(p_i,p_j),
$$
or equivalently,  the matrix $\bigl(f(t_{ij})\bigr)\succeq0$, where  the sign $\succeq 0$ stands for: ``is  positive semidefinite".

Let ${\Bbb S}\sp {n-1}$ denote the unit sphere in ${\Bbb R}\sp n$, and let $\varphi_{ij}$ denote the angular distance between points $p_i, p_j.$  Schoenberg \cite{Scho} proved that:\\
\centerline{\em $f(\cos\varphi)$ is p.d. in ${\Bbb S}\sp {n-1}$ if and only if
$\; f(t)=\sum_{k=0}^\infty{f_kG_k^{(n)}(t)}$ with all $f_k\ge 0$.}
\\
 Here $G_k^{(n)}(t)$
are the Gegenbauer   polynomials.

Schoenberg's theorem has been generalized by Bochner \cite{Boc} to more general spaces. Namely, the following fact holds: {\em $f$ is p.d. in a 2-point-homogenous space $M$ if and only if $f(t)$ is a nonnegative linear combination   of the zonal spherical functions  $\Phi_k(t)$} (see details in \cite[Th. 2]{Kab}, \cite[Chapter 9]{CS}).

Note that
the Bochner - Schoenberg theorem is widely used in coding theory and discrete geometry for finding bounds for error-correcting codes, constant weight codes, spherical codes, sphere packings and other packing problems in 2-point-homogeneous spaces (see \cite{CS,  Kab, Mus2,Mus3, Mus4, PZ} and many others).

Consider the following problem: Let {\em  $Q=\{q_1,\ldots,q_m\}$ be a set of points in  $M$.   To describe the class  of  continuous functions $F(t,u,v)$ in $2m+1$ variables with $ t\in {\Bbb R},\; u,v\in {\Bbb R}\sp m, \; F(t,u,v)=F(t,v,u)$ such that for arbitrary points $p_1,\ldots,p_r$ in $M$ the matrix
$
\bigl(F(t_{ij},u_i,u_j)\bigr)\succeq 0, \; \mbox { where } \; t_{ij}=\tau(p_i,p_j), \;
u_i=(\tau(p_i,q_1),\ldots,\tau(p_i,q_m)).\; $} 

Denote this class by $\psd(M,Q)$. If
$Q=\emptyset$,  then $\psd(M,Q)$ is the class of p.d. functions in $M$.



Recently, Schrijver \cite{Schr}  improved some upper bounds on binary codes using  semidefinite programming.
Schrijver's method has been adapted
 for some non-binary codes (Gijswijt,  Schrijver, and Tanaka \cite{GST}), and  for spherical codes  (Bachoc and Vallentin \cite{BV, BV2, BV3}).   

In fact, by using the stabilizer subgroup of the isometry group  this method  derives new positive semidefinite constraints which are stronger than constraints given by the Bochner - Schoenberg theorem (see \cite{BV}). In other words, it shows that constraints for $\{t_{ij}\}$ and $\{u_{i}\}$ given by
$\psd(M,\{q_1\})$ are stronger than constraints in $\psd(M,\emptyset)$ for spherical, Hamming and some other spaces.

Clearly, if $F_1, F_2\in\psd(M,Q)$, then for any functions $\alpha_1, \alpha_2$ in $m$ variables we have
$$F(t,u,v)=\alpha_1(u)\alpha_1(v)F_1(t,u,v)+\alpha_2(u)\alpha_2(v)F_2(t,u,v)\in\psd(M,Q). \eqno (1.1)$$

In this paper we consider the class $\psd({\Bbb S}\sp {n-1},Q)$. Namely, we show that this class can be generated by  $(1.1)$ using as a basis the polynomials
$G_k^{(n-m)}$,  where $m=\card(Q)\le n-2.$  Actually, the class $\psd({\Bbb S}\sp {n-1},Q)$ becomes ``more rich" whenever $m$ is increasing.


The paper is organized as follows: Section 2 introduces multivariate Gegenbauer polynomials and its basic properties. Section 3 extends the Schoenberg theorem. Section 4 extends the Schoenberg theorem to Euclidean spaces.  Section 5  considers the positive semidefinite constraints for the distance distribution given by Theorem 3.1.
Section 6 discusses applications of the extended Bochner - Schoenberg theorem to spherical codes.

\section{Multivariate Gegenbauer polynomials}
In this section we introduce polynomials $G_k^{(n,m)}(t,{\bf u},{\bf v})$ in $2m+1$ variables with $0\le m\le n-2,\; t\in {\Bbb R},\;$
${\bf u}=(u_1,\ldots,u_m),$
${\bf v}=(v_1,\ldots,v_m)\in {\Bbb R}\sp {m}$.

\medskip

\noindent{\bf 2.1. Gegenbauer polynomials.}
There are many ways to define Gegenbauer (or ultraspherical) polynomials
$G_k^{(n)}(t)$ (see \cite{CS,Erd,PZ,Scho}).  $G_k^{(n)}$ are a special case of  Jacobi polynomials
$P_k^{(\alpha,\beta)}$ with $\alpha=\beta=(n-3)/2$ and the normalization by $G_k^{(n)}(1)=1$.
Also
$G_k^{(n)}$ can be defined by the recurrence formula:
$$ G_0^{(n)}=1,\; G_1^{(n)}=t,\; \ldots,\; G_k^{(n)}=\frac {(2k+n-4)\,t\,G_{k-1}^{(n)}-(k-1)\,G_{k-2}^{(n)}} {k+n-3}.$$

Note that for any even $k\ge0$ (resp. odd)  $G_k^{(n)}(t)$ is even (resp. odd). Therefore,
$G_{2k}$ and $G_{2\ell+1}$ are orthogonal on $[-1,1]$.
Moreover, all  polynomials $G_k^{(n)}$ are orthogonal on $[-1,1]$ with respect to the weight
 function $(1-t^2)^{(n-3)/2}$:
$$
\int_{-1}^1G_k^{(n)}(t)\,G_\ell^{(n)}(t)\,(1-t^2)^{(n-3)/2}\,dt=0, \quad k\ne\ell. \eqno (2.1)$$

Recall the { addition theorem for Gegenbauer polynomials}:
$$
G_k^{(n)}(\cos{\theta_1}\cos{\theta_2}+\sin{\theta_1}\sin{\theta_2}\cos{\varphi})$$
$$=\sum\limits_{s=0}^k c_{nks}\,G_{k-s}^{(n+2s)}(\cos{\theta_1})\,G_{k-s}^{(n+2s)}(\cos{\theta_2})\,(\sin{\theta_1})^s\,(\sin{\theta_2})^s\,G_s^{(n-1)}(\cos{\varphi}),
$$
where $c_{nks}$ are  {\em positive} coefficients whose values of  no concern here  (see \cite{Car,Erd}).

\medskip

\noindent{\bf 2.2. Basic definitions.}
Let us denote by $\langle {\bf u},{\bf v}\rangle$  the inner product of vectors
${\bf u},{\bf v}\in{\Bbb R}\sp m$, and  by $|{\bf v}|$ we denote  the norm of ${\bf v}$ (i.e. $|{\bf v}|^2=\langle {\bf v},{\bf v}\rangle$).

\begin{defn}
Let $0\le m\le n-2, \; t\in {\Bbb R},\; {\bf u},{\bf v}\in{\Bbb R}\sp m$ for $m>0$, and
${\bf u}={\bf v}=0$ for $m=0.$
 Then the following polynomial  in $2m+1$ variables of degree $k$ in $t$ is well defined:
$$
G_k^{(n,m)}(t,{\bf u},{\bf v}):=(1-|{\bf u}|^2)^{k/2}\,(1-|{\bf v}|^2)^{k/2}\,
G_k^{(n-m)}\left(\frac{t-\langle{\bf u},{\bf v}\rangle}
{\sqrt{(1-|{\bf u}|^2)(1-|{\bf v}|^2)}}\right).
$$

\end{defn}

Note that by definition we have: $ G_k^{(n,0)}(t)=G_k^{(n,0)}(t,0,0)=G_k^{(n)}(t).$

Let ${\bf u}=(u_1,\ldots,u_{n})\in {\Bbb R}\sp {n}.$ For $0<m\le n$ we denote by ${\bf u}^{(m)}$ the vector
$(u_1,\ldots,u_m)$ in ${\Bbb R}\sp {m}$. Let ${\bf u}^{(0)}=0$.
 For $0\le m\le n-2$, ${\bf u},{\bf v}\in {\Bbb R}\sp {n}$ put
$$G_k^{(n,m)}(t,{\bf u},{\bf v}):=G_k^{(n,m)}(t,{\bf u}^{(m)},{\bf v}^{(m)}).$$

\begin{defn}
We define a matrix $Z_d^m({\bf u},{\bf v})$ 
of size ${m+d\choose m}\times{m+d\choose m}$ for
${\bf u}\in{\Bbb R}\sp m,\; {\bf v}\in{\Bbb R}\sp m$
 by
$\; Z_d^m({\bf u},{\bf v}):=(z_d^m({\bf u}))^T\,z_d^m({\bf v}),\; $ where
$$z_d^m({\bf x})=z_d^m(x_1,\ldots,x_m):=(1,x_1,\ldots,x_m,x_1^2,x_1x_2,\ldots,x_{m-1}x_m,x_m^2,\ldots,x_m^d)$$ is the vector of monomials.
\end{defn}

Denote   the inner product     of  matrices $A=\bigl(a_{ij}\bigr),\; B=\bigl(b_{ij}\bigr)$ both of size $d\times d$  by
$\langle A,B\rangle$, i.e.
$\langle A,B\rangle=\tr(AB)=\sum\limits_{i,j=1}^d {a_{ij}b_{ij}}$.

\begin{defn} Let  $f({\bf u},{\bf v})$ be a  symmetric polynomial in ${\bf u},{\bf v}\in{\Bbb R}\sp m$, i.e. $f({\bf u},{\bf v})=f({\bf v},{\bf u})$.
We say that $f({\bf u},{\bf v})$ is positive semidefinite and write $f\succeq 0$ if there are a symmetric matrix $H\succeq0$ and $d\ge0$ so that   $f({\bf u},{\bf v})=\left\langle H,Z_{d}^{m}({\bf u},{\bf v})\right\rangle$.

Let polynomials $f_r\succeq0$ for all $r, \; f_r\to f$ as $r\to\infty$, where a function $f({\bf u},{\bf v})$ is
symmetric in ${\bf u},{\bf v}$  and continuous on $\{|{\bf u}|\le 1, |{\bf v}|\le 1\}$. Then we say that $f({\bf u},{\bf v})$ is positive semidefinite and write $f\succeq 0$.
\end{defn}

\noindent{\bf 2.3. Orthogonality.}
Let $t\in {\Bbb R},\;$
${\bf u}=(u_1,\ldots,u_m),$
${\bf v}=(v_1,\ldots,v_m)\in {\Bbb R}\sp {m}$.
Denote
$$
Q_m(t,{\bf u},{\bf v}):=
\left(
\begin{array}{ccccc}
1 & \ldots & 0  &  u_{1} &  v_{1}\\
\vdots & \vdots
&\ddots & \vdots\\
0 & \ldots & 1  &  u_{m} &  v_{m}\\
 u_{1} &\ldots &  u_{m} &  1 &  t\\
 v_{1}
&\ldots &  v_{m} &  t &  1
\end{array}
\right).
$$
Let
$$
D_m:=\{(t,{\bf u},{\bf v}) \in {\Bbb R}\sp {2m+1}: \; Q_m(t,{\bf u},{\bf v})\succeq 0\}.
$$
Since
$$
\det(Q_m(t,{\bf u},{\bf v}))=(1-|{\bf u}|^2)(1-|{\bf v}|^2)-(t-\langle {\bf u},{\bf v}\rangle)^2,
$$
we have $\; (t,{\bf u},{\bf v})\in D_m  \; $ if and only if
$$
(t-\langle {\bf u},{\bf v}\rangle)^2\le(1-|{\bf u}|^2)(1-|{\bf v}|^2),\quad |{\bf u}|\le 1.
$$

Denote $(\det(Q_m(t,{\bf u},{\bf v})))^{(n-m-3)/2} \; $ by $\; \rho_{n,m}(t,{\bf u},{\bf v})$, i.e.
$$
\rho_{n,m}(t,{\bf u},{\bf v}):=
\bigl((1-|{\bf u}|^2)(1-|{\bf v}|^2)-(t-\langle {\bf u},{\bf v}\rangle)^2\bigr)^{(n-m-3)/2}.
$$

\begin{theorem} Let $\; 0\le m\le n-2. \; $ Let
$q({\bf u},{\bf v})$  be any continuous function on $\{({\bf u},{\bf v})\in{\Bbb R}\sp {2m}: |{\bf u}|\le 1, |{\bf v}|\le 1\}.\; $
If $\; k\ne \ell, \; $ then
$$
\int_{D_m}{G_{k}^{(n,m)}(t,{\bf u},{\bf v})\,G_{\ell}^{(n,m)}(t,{\bf u},{\bf v})\,
q({\bf u},{\bf v})\,\rho_{n,m}(t,{\bf u},{\bf v})\,dt\,d{\bf u}\,d{\bf v}}=0.
$$
\end{theorem}
\begin{proof} Let $$t=s\,({(1-|{\bf u}|^2)(1-|{\bf v}|^2)})^{1/2}+
\langle {\bf u},{\bf v}\rangle. $$
In the variables $s, {\bf u}, {\bf v}$ we have
$$D_m=\{(s,{\bf u},{\bf v}): -1\le s \le 1,\; |{\bf u}|\le 1,\; |{\bf v}|\le 1\}.$$ The Jacobian of this change of variables is
$({(1-|{\bf u}|^2)(1-|{\bf v}|^2)})^{1/2}$.
Then
$$
I=\int_{D_m}{G_{k}^{(n,m)}(t,{\bf u},{\bf v})\,G_{\ell}^{(n,m)}(t,{\bf u},{\bf v})\,
q({\bf u},{\bf v})\,\rho_{n,m}(t,{\bf u},{\bf v})\,dt\,d{\bf u}\,d{\bf v}} $$ $$
=
I_0\,\int_{-1}^1G_k^{(n-m)}(s)\,G_\ell^{(n-m)}(s)\,(1-s^2)^{(n-m-3)/2}\,ds,
$$
$$
I_0=\int_{|{\bf u}|\le 1}\int_{|{\bf v}|\le 1}\bigl((1-|{\bf u}|^2)(1-|{\bf v}|^2)\bigr)^{(n-m+k-2)/2}\,q({\bf u},{\bf v})\,
d{\bf u}\,d{\bf v}.
$$
Thus, $(2.1)$ yields $I=0.$
\end{proof}

Let $e_1, \ldots, e_n$ be an orthonormal basis of ${\Bbb R}\sp n$, and let $(x_1,\ldots,x_n)$ be the coordinates expression  of a point $x\in{\Bbb R}\sp n$ in this basis.  Then $x^{(m)}$ as well as
$G_k^{(n,m)}(\langle x,y\rangle,x,y)=G_k^{(n,m)}(\langle x,y\rangle,x^{(m)},y^{(m)})$
for $x, y\in{\Bbb R}\sp n$ are well defined.

\begin{lemma}
For any continuous $F$ on $D_m$ and $\; 0\le m\le n-2$ we have
$$
\int_{{({\Bbb S} \sp {n-1})^2}}
{F\bigl(\langle x,y\rangle,x^{(m)},y^{(m)}\bigr)\,d\omega_n(x)\,d\omega_n(y)}$$
$$
=\omega_{n-m}\,\omega_{n-m-1}\int_{D_m}{F(t,{\bf u},{\bf v})\,r_{m}({\bf u})\,r_{m}({\bf v})\,\rho_{n,m}(t,{\bf u},{\bf v})\,dt\,d{\bf u}\,d{\bf v}},
$$
where $\; r_{m}({\bf u})=\prod_{i=1}^{m-1}({1-|{\bf u}^{(i)}|^2})^{1/2}$ for $m>1,\;
r_0=r_1=1$, and
 $\omega_n$ is the surface area of ${\Bbb S}\sp {n-1}$ for the standard measure $d\omega_n$.
\end{lemma}
%
%
\begin{proof}  Let  $x,y\in{\Bbb S}\sp {n-1}$.
Let $t=\langle x,y\rangle,\;  {\bf u}=x^{(m)},\; {\bf v}=y^{(m)}$. Let
$$
a=\frac{x-{\bf u}}{h({\bf u})},\quad b=\frac{y-{\bf v}}{h({\bf v})},\quad
s=\langle a,b\rangle, \; \mbox{ where } \;
h({\bf u})=\sqrt{1-|{\bf u}|^2}.
$$
Then
$$
{F\bigl(\langle x,y\rangle,x^{(m)},y^{(m)}\bigr)}=
F(h({\bf u})\,h({\bf v})\,\langle a,b\rangle +\langle {\bf u},{\bf v}\rangle,{\bf u},{\bf v})=
F(t,{\bf u},{\bf v}).
$$
Therefore, we can change the variables
$
(x,y) \to (a, b, {\bf u}, {\bf v}) \to (s,{\bf u}, {\bf v})  \to (t, {\bf u}, {\bf v}).
$

 We obviously  have $|{\bf u}|\le 1, \;|{\bf v}|\le 1,\, -1\le s\le1$. Clearly, $|a|=|b|=1$. Then
$(a,b)\in({\Bbb S}\sp {n-m-1})^2. \; $ Since $|t-\langle {\bf u},{\bf v}\rangle|\le
h({\bf u})\,h({\bf v})$ we have
$(t,{\bf u},{\bf v})\in D_m$.

 Consider
$(x,y) \to (a, b, {\bf u}, {\bf v})$.   The Jacobian is defined by
$$
d\omega_n(x)=k_{n,m}({\bf u})\,d{\bf u}\,d\omega_{n-m}(a), $$ $$
k_{n,m}({\bf u})=({1-|{\bf u}|^2})^{(n-m-2)/2}\,
\prod_{i=0}^{m-1}({1-|{\bf u}^{(i)}|^2})^{1/2}.
$$
Combining the volume of $\{(a,b)\in({\Bbb S}\sp {n-m-1})^2: \langle a,b\rangle=s\}$ and the Jacobian of $(a, b, {\bf u}, {\bf v}) \to (s,{\bf u}, {\bf v})$ we get:
$\; \omega_{n-m}\,\omega_{n-m-1}\,
\rho_{n,m}(s,0,0).$
The Jacobian of the changing $(s,{\bf u},{\bf v})  \to (t,{\bf u},{\bf v})$ is
$(h({\bf u})\,h({\bf v}))^{-1}$.
That completes the proof.
\end{proof}


Combining Lemma 2.1 and Theorem 2.1 we get:
\begin{cor} Let $\; 0\le m\le n-2.$ Let
$f({\bf u},{\bf v})$  be any continuous function  on $\{({\bf u},{\bf v})\in{\Bbb R}\sp {2m}: |{\bf u}|\le 1, |{\bf v}|\le 1\}.\; $   If $\; k\ne \ell, \; $ then
$$
\int_{{({\Bbb S}\sp {n-1})^2}}
{G_{k}^{(n,m)}(\langle x,y\rangle,x,y)\,G_{\ell}^{(n,m)}(\langle x,y\rangle,x,y)}\,
f\bigl(x^{(m)},y^{(m)}\bigr)\,d\omega_n(x)\,d\omega_n(y)
=0.
$$
\end{cor}

\medskip

\noindent{\bf 2.4. The addition theorem.}
\begin{theorem}[The addition theorem] Let $1\le m \le n-2$. Then
$$
G_k^{(n,m-1)}(t,{\bf u},{\bf v})=\sum\limits_{s=0}^k{C_{k-s}^{n,m}({\bf u})\,C_{k-s}^{n,m}({\bf v})\,G_s^{(n,m)}(t,{\bf u},{\bf v})},
$$
where $C_{d}^{n,m}({\bf u})=C_{d}^{n,m}(u_1,\ldots,u_m)$ is a polynomial in $u_1,\ldots,u_m$ of degree $d$.
\end{theorem}
\begin{proof}
Suppose  $t,\; {\bf u}=(u_1,\ldots,u_m),\; {\bf v}=(v_1,\ldots,v_m)$ are such that
$$|{\bf u}|< 1,  \quad |{\bf v}|< 1,\quad
 (t-\langle {\bf u},{\bf v}\rangle)^2 < {(1-|{\bf u}|^2)(1-|{\bf v}|^2)}.$$
 Then  $\; \theta_1, \theta_2, \varphi\in (0,\pi) $ 
 are uniquely defined by the following equations:
$$
\cos{\theta_1}=\frac{u_{m}}{\sqrt{1-u_1^2-\ldots-u_{m-1}^2}},\quad
\cos{\theta_2}=\frac{v_{m}}{\sqrt{1-v_1^2-\ldots-v_{m-1}^2}},
$$
$$
\cos{\varphi}=\frac{t-\langle {\bf u},{\bf v}\rangle}{\sqrt{(1-|{\bf u}|^2)(1-|{\bf v}|^2)}}.
$$

The addition theorem for Gegenbauer polynomials yields:
$$
G_k^{(n,m-1)}(t,{\bf u},{\bf v})=\sum\limits_{s=0}^k{C_{k-s}^{n,m}({\bf u})\,C_{k-s}^{n,m}({\bf v})\,G_s^{(n,m)}(t,{\bf u},{\bf v}) },\eqno (2.2)
$$
$$
C_{k-s}^{n,m}({\bf u}):=\sqrt{c_{nks}}\,w^{k-s}\,G_{k-s}^{(n+2s)}\left(\frac{u_{m}}{w}\right),
\quad w=\sqrt{1-u_1^2-\ldots-u_{m-1}^2}.
$$
It's easy to see that $C_{k-s}^{n,m}(u_1,\ldots,u_m)$ is a polynomial of degree $k-s$. Thus, $(2.2)$ holds for all $t, {\bf u}, {\bf v}.$
\end{proof}





Let $b_s({\bf u},{\bf v})=C_{k-s}^{n,m}({\bf u})\,C_{k-s}^{n,m}({\bf v})$. Clearly, $b_s\succeq 0$. Then Theorem 2.2 yields
%

\begin{cor} Let $0\le m \le \ell \le n-2.$ Let ${\bf u},{\bf v}\in{\Bbb R}\sp {n}$.   Then
$$
G_k^{(n,m)}(t,{\bf u},{\bf v})
=\sum\limits_{s=0}^k{f_s\bigl({\bf u}^{(\ell)},{\bf v}^{(\ell)}\bigr)\,G_s^{(n,\ell)}(t,{\bf u},{\bf v})}
$$
with $f_s\succeq0$ for all $\; 0\le s\le k.$ 
\end{cor}

\section{An extension of the Schoenberg theorem}
{\bf 3.1. Schoenberg's theorem.}
Let  $p_1, \ldots, p_r$ be  points in ${\Bbb S}\sp {n-1}$, and let $a_1,\ldots, a_r$ be any real numbers. Then
$$ 0\le {\left|a_1p_1+\ldots +a_rp_r\right|}^2 = \sum\limits_{i,j}
{\langle p_i,p_j\rangle}a_ia_j,$$ or equivalently
the Gram matrix $\Big(\langle p_i,p_j\rangle\Big)$ is  positive semidefinite.

Schoenberg \cite{Scho} extended this property to Gegenbauer polynomials $G_k^{(n)}$.
He proved that {\em for any finite $X=\{p_1,\ldots, p_r\}\subset {\Bbb S}\sp {n-1}$
 the matrix  $\Big(G_k^{(n)}(\langle p_i,p_j\rangle)\Big)$  is positive semidefinite.}

{
Schoenberg proved also that the converse holds: {\em if $f(t)$ is a real polynomial and for any finite $X\subset{\Bbb S}\sp {n-1}$ the matrix $\big(f(\langle p_i,p_j\rangle)\big)\succeq 0$, then $f(t)$ is a linear combination of $G_k^{(n)}(t)$ with nonnegative coefficients.}}

\medskip


\noindent{\bf 3.2. An extension of the direct Schoenberg theorem.}
\begin{theorem}  Let $e_1, \ldots, e_n$ be an orthonormal basis of ${\Bbb R}\sp n$, and let $p_1,\ldots, p_r$  be points in ${\Bbb S}\sp {n-1}$.
Then for any $k\ge0$ and  $0\le m \le n-2$  the matrix
$\Bigl(G_k^{(n,m)}(\langle p_i,p_j\rangle,p_i,p_j)\Bigr)$ is positive semidefinite.
\end{theorem}

\begin{proof}
Actually, this theorem is a simple consequence of the
Schoenberg theorem.
Indeed, let $ {\bf v}_i=\langle p_i,e_1\rangle\,e_1+\ldots+\langle p_i,e_m\rangle\,e_m$,
and let ${\bf x}_i=p_i-{\bf v}_i$. Then ${\bf x}_i$ is a vector in ${\Bbb R}\sp {n-m}$ with the basis
$e_{m+1}, \ldots, e_n$. Note that  $|{\bf v}_i|^2+|{\bf x}_i|^2=1$.

Let  ${\bf y}_i={\bf x}_i/|{\bf x}_i|$  for $|{\bf x}_i|>0$.
In the case $|{\bf x}_i|=0$ put ${\bf y}_i=e_n.$
Then ${\bf y}_i\in{\Bbb S}\sp {n-m-1}$.

Recall the Schur theorem: {\em If $A=\bigl(a_{ij}\bigr)\succeq 0, \; B=\bigl(b_{ij}\bigr)\succeq 0$, then $C=\bigl(a_{ij}b_{ij}\bigr)\succeq 0. \; $}
Let
$ h_i=|{\bf x}_i|^k=(1-|{\bf v}_i|^2)^{k/2},  \;
{\bf h}=\left(h_1,\ldots,h_r\right), \; A={\bf h}^T{\bf h}.
$
Clearly, $A\succeq0.$
The Schoenberg theorem yields:
 $B=\Bigl({G_k^{(n-m)}(\langle {\bf y}_i,{\bf y}_j\rangle)}\Bigr)\succeq 0.$
It is easy to see that
$G_k^{(n,m)}(t_{ij},{\bf v}_i,{\bf v}_j)=a_{ij}b_{ij},$
where $t_{ij}:=\langle p_i,p_j\rangle$.
Thus, $C=\Bigl(G_k^{(n,m)}(t_{ij},{\bf v}_i,{\bf v}_j)\Bigr)\succeq 0$.

\medskip

Using the addition theorem it is not hard to give a direct proof. Note that $G_k^{(2)}$ is the Chebyshev polynomial of the first kind, i.e. $G_k^{(2)}(\cos\theta)=\cos(k\theta).$ It is easy to see that for any $\varphi_1,\ldots,\varphi_r$ the matrix
$\Bigl(\cos(\varphi_i-\varphi_j)\Bigr)\succeq0.$ From this follows that Schoenberg's theorem holds for $n=2.$ Therefore, we have proved the theorem for $m=n-2.$  Put $\ell=n-2$ in Corollary 2.2. That yields the theorem for all $0\le m\le n-2.$
 \end{proof}

\noindent{\bf Remark 3.1.} We see that the second proof (as well as Schoenberg's original proof in \cite{Scho}) is based on the addition theorem for Gegenbauer polynomials. There exists another proof of Schoenberg's theorem which is using the {\em addition theorem for spherical harmonics} (see, for instance, \cite{PZ}). It is possible (see \cite{BV} for the case $m=1$) to derive Theorem 3.1 from this theorem. However, this proof looks more complicated than our proof.

\medskip




\noindent{\bf Remark 3.2.} Bachoc and Vallentin \cite{BV} derived new upper bounds for spherical codes based on  positive semidefinite constraints that are given in \cite[Corollary 3.4]{BV}.  In fact, this corollary
easily follows from Theorem 3.1 with $m=1.$

Indeed, let us fix some $\ell, \; 1\le\ell\le r$, and take an orthonormal basis $e_1, \ldots, e_n$ of ${\Bbb R}\sp n$ with $e_1=p_\ell$. Then  for $m=1$ we have ${\bf v}_i=t_{i\ell}$. Let
$$\left(A_\ell\right)_{ij}=G_k^{(n,1)}(t_{ij},t_{i\ell},t_{j\ell}),
\quad  0\le i\le r, \; 1\le j\le r,$$
$$Y_\ell=WA_\ell W^T, \quad
(W)_{ij}=\lambda_iG_i^{(n+2k)}(t_{j\ell}),\; 0\le i\le d, \; 1\le j\le r.
$$
Theorem 3.1 yields: $A_\ell \succeq 0$.
Then $Y_\ell\succeq 0$ and $Y_1+\ldots+Y_r\succeq 0$. These constraints are equivalent to   the constraints in \cite[Corollary 3.4]{BV}.

\begin{cor}  Let $0\le m\le n-2, \; d\ge0$.  Let $e_1, \ldots, e_n$ be an orthonormal basis of ${\Bbb R}\sp n$.
Let a polynomial $F(t,{\bf u},{\bf v})$ can be represented in the form
$$
F(t,{\bf u},{\bf v})=\sum\limits_{k=0}^d{f_k({\bf u},{\bf v})\,G_{k}^{(n,m)}(t,{\bf u},{\bf v})},
$$
where $f_k\succeq0$  for all $\; 0\le k\le d.$
Then for any points  $p_1,\ldots, p_r$  in ${\Bbb S}\sp {n-1}$
the matrix
$\Bigl(F\bigl(\langle p_i,p_j\rangle,p_i^{(m)},p_j^{(m)}\bigr)\Bigr)$ is positive semidefinite.
In other words, $F\in\psd({\Bbb S}\sp {n-1},\{e_1,\ldots,e_m\})$.
\end{cor}
\begin{proof} 
%
Let $f({\bf u},{\bf v})\succeq0.$ Then $f({\bf u},{\bf v})=\left\langle H,Z_{d}^{m}({\bf u},{\bf v})\right\rangle$ with
  $H\succeq 0$. From an eigenvalue factorization of $H$ it follows that there exist polynomials $h_i({\bf u})$ such that $f({\bf u},{\bf v})=\sum_i{h_i({\bf u})\,h_i({\bf v})}$.
Therefore,  $\Bigl(f\bigl(p_i^{(m)},p_j^{(m)}\bigr)\Bigr)\succeq0$.

We have
$A_k=\Bigl(f_k\bigl(p_i^{(m)},p_j^{(m)}\bigr)\Bigr)\succeq0.$
Theorem 3.1 yields
$$B_k=\Bigl(G_k^{(n,m)}\bigl(\langle p_i,p_j\rangle,p_i^{(m)},p_j^{(m)}\bigr)\Bigr)\succeq0.$$

Let $(C_k)_{ij}=(A_k)_{ij}\,(B_k)_{ij}$.  The Schur theorem implies: $C_k\succeq0.\; $ Thus, $\Bigl(F\bigl(\langle p_i,p_j\rangle,p_i^{(m)},p_j^{(m)}\bigr)\Bigr)=\sum\limits_{k=0}^d {C_k}\succeq0.$
\end{proof}

\noindent{\bf 3.3. An extension of the converse Schoenberg theorem.}
\begin{theorem}  Let $0\le m\le n-2$.  Let $e_1, \ldots, e_n\in {\Bbb R}\sp n$ be an orthonormal basis.
Let $F(t,{\bf u},{\bf v})$ be a polynomial in $t$ and a symmetric polynomial in
${\bf u},{\bf v}\in {\Bbb R}\sp m$. 
Suppose
$F\in\psd({\Bbb S}\sp {n-1},\{e_1,\ldots,e_m\})$. Then
$$F(t,{\bf u},{\bf v})=\sum\limits_{k=0}^d{f_k({\bf u},{\bf v})\,G_{k}^{(n,m)}(t,{\bf u},{\bf v})}
\eqno (3.1)
$$
with $f_k\succeq0$ for all $\; 0\le k\le d=\deg_t(F).$
\end{theorem}
\begin{proof} The polynomial $G_k^{(n,m)}(t,{\bf u},{\bf v})$ has degree $k$ in the variable $t$, so that $F$ has a unique expression in the form $(3.1)$,
where $f_k({\bf u},{\bf v})$ is a symmetric polynomial in ${\bf u},{\bf v}$ of degree $d_k$ in ${\bf u}$. Then $f_k({\bf u},{\bf v})=\left\langle F_k,Z_{d_k}^{m}({\bf u},{\bf v})\right\rangle$, where $F_k$ is a symmetric matrix.
Now we prove that $F_k\succeq0.$

Let  $f(t,{\bf u},{\bf v}), \; g(t,{\bf u},{\bf v})$ be continuous functions on $D_m$. Denote
$$
\Bigl\{f,g\Bigr\}:=
\int_{{({\Bbb S}\sp {n-1})^2}}
f(\langle x,y\rangle,x^{(m)},y^{(m)})\,g(\langle x,y\rangle,x^{(m)},y^{(m)})\,
d\omega_n(x)\,d\omega_n(y).
$$
For any continuous function $h({\bf u},{\bf v})$ Corollary 2.1 yields:
$$
\Bigl\{F,h\,G_k^{(n,m)}\Bigr\}=\Bigl\{f_k\,G_k^{(n,m)},h\,G_k^{(n,m)}\Bigr\}=
\left\{f_k\,h,\bigl(G_k^{(n,m)}\bigr)^2\right\}. \eqno (3.2)
$$

Let $f({\bf u},{\bf v}), \; g({\bf u},{\bf v})$ be continuous functions on
$\{|{\bf u}|\le1,\; |{\bf v} |\le1\}$. Denote
$$
[f,g]:=\left\{fg,\bigl(G_k^{(n,m)}\bigr)^2\right\}.
$$
Clearly, $[\cdot,\cdot]$ is an inner product.  Let $\alpha_1, \alpha_2,\ldots$ be an orthonormal basis for $[\cdot,\cdot]$ in the space of real polynomials  ${\bf {\Bbb R}[u,v}]$. We observe that $F_k$ defines a quadratic form on  ${\bf {\Bbb R}[u,v}]$. Let us denote by  $\tilde F_k$  the matrix expression of this quadratic form in the basis  $\{\alpha_i\}.$ Then there exists a matrix $A$ such that $F_k=A^T\tilde F_kA$. Therefore, $\tilde F_k\succeq 0$ yields $F_k\succeq 0$.

Let $h({\bf u},{\bf v})=\left\langle H,Z_{d_k}^{m}({\bf u},{\bf v})\right\rangle$.  It's easy to see that $[f_k,h]=\langle\tilde F_k,\tilde H\rangle$.

Note that the sum of all entries of a positive semidefinite matrix is nonnegative. This implies: if $f\in\psd({\Bbb S}\sp {n-1},E_m),\; E_m:=\{e_1,\ldots,e_m\}$, then $\bigl\{f,1\bigr\}\ge 0$. Moreover,
the Schur theorem yields: if $f$ and $g\in\psd({\Bbb S}\sp {n-1},E_m)$, then $fg\in\psd({\Bbb S}\sp {n-1},E_m)$. Therefore,
$\bigl\{f,g\bigr\}\ge0.$ Thus, for any  $f, g\in\psd({\Bbb S}\sp {n-1},E_m)$ we have $[f,g]\ge0.$

Let
$H\succeq0.$ Then
$h({\bf u},{\bf v})\in\psd({\Bbb S}\sp {n-1},E_m)$.
We have $F\in\psd({\Bbb S}\sp {n-1},E_m)$. Then $(3.2)$ yields $[f_k,h]=\langle\tilde F_k,\tilde H\rangle\ge0$
for any $\tilde H\succeq0$. Thus, $\tilde F_k\succeq 0$ and $F_k\succeq 0$.
\end{proof}

\noindent{\bf Remark 3.3.} In \cite[Proposition 4.12]{BV2} a statement that is equivalent to Theorem 3.2 was proven  with $m=1$.   However, in \cite[Proposition 4.12]{BV2} the formula $(3.1)$ is written in terms of an orthogonal  basis for $\bigl\{\cdot,\cdot\bigr\}$.

\medskip

\noindent{\bf  3.4.  The class $\psd({\Bbb S}\sp {n-1},Q)$.} If $f_r\in\psd(M,Q),\; f_r \to f$ as
$r\to\infty$, and $f(t,{\bf u},{\bf v})$ is continuous, then also $f\in\psd(M,Q)$. Therefore, Corollary 3.1 and Theorem 3.2 imply
\begin{theorem}  Let $0\le m\le n-2$.  Let $e_1, \ldots, e_n$ be an orthonormal basis of ${\Bbb R}\sp n$.
Then
$F\in\psd({\Bbb S}\sp {n-1},\{e_1,\ldots,e_m\})$ if and only if
$$
F(t,{\bf u},{\bf v})=\sum\limits_{k=0}^\infty{f_k({\bf u},{\bf v})\,G_{k}^{(n,m)}(t,{\bf u},{\bf v})},
$$
where  for any $k\ge 0$ a function $f_k({\bf u},{\bf v})$ is positive semidefinite.
\end{theorem}

Now we consider $\psd({\Bbb S}\sp {n-1},Q)$ for any $Q=\{q_1,\ldots,q_m\}\subset{\Bbb S}\sp {n-1}$ with  $0\le m\le n-2$. Theorem 3.3 yields
\begin{cor}
Let $0\le m\le n-2$. Let $Q=\{q_1,\ldots,q_m\}\subset{\Bbb S}\sp {n-1}$ with $\rk(Q)=m$.
Let $e_1, \ldots, e_m$ be an orthonormal basis
  of the linear space
 with the basis $q_1,\ldots,q_m$, and
let $L_Q$ denotes the linear transformation of coordinates.
Then
$F\in\psd({\Bbb S}\sp {n-1},Q)$ if and only if
$$
F(t,{\bf u},{\bf v})=\sum\limits_{k=0}^\infty{f_k({\bf u},{\bf v})\,G_{k}^{(n,m)}(t,L_Q({\bf u}),L_Q({\bf v}))},
$$
where   $f_k({\bf u},{\bf v})\succeq0$ for all $k\ge 0$.
\end{cor}

\noindent{\bf Remark 3.4.} It is not hard to describe $\psd({\Bbb S}\sp {n-1},Q)$ with $\rk(Q)\ge n-1.$

First, consider the case $Q=\{e_1,\ldots,e_{n-1}\}$, where $\{e_i\}$ is an orthonormal basis.
Let  $p_1,\ldots,p_r\in{\Bbb S}\sp {n-1},$   ${\bf u}_i=\langle p_i,e_1\rangle\,e_1+\ldots+\langle p_i,e_{n-1}\rangle\,e_{n-1}$, and $t_{ij}=\langle p_i,p_j\rangle.$
Obviously, $t_{ij}-\langle {\bf u}_i,{\bf u}_j\rangle=w_iw_j, \;
 w_i=\langle p_i,e_n\rangle$. Then
$$\Bigl(t_{ij}-\langle {\bf u}_i,{\bf u}_j\rangle\Bigr)={\bf w}^T{\bf w}\succeq0, \quad {\bf w}=(w_1,\ldots,w_r). \eqno (3.3) $$
Since $w_i^2=1-|{\bf u}_i|^2$, we have for all $i=1,\ldots,r,\; j=1,\ldots,r$:
$$(t_{ij}-\langle {\bf u}_i,{\bf u}_j\rangle)^2=(1-|{\bf u}_i|^2)(1-|{\bf u}_j|^2).\eqno (3.4)$$
Let $
H_2(t,{\bf u},{\bf v}):=(t-\langle {\bf u},{\bf v}\rangle)^2-(1-|{\bf u}|^2)(1-|{\bf v}|^2).
$ Since $H_2(t_{ij},{\bf u}_i,{\bf u}_j)=0$, we have
 $\; F\in\psd({\Bbb S}\sp {n-1},\{e_1,\ldots,e_{n-1}\})\; $ if and only if
$$
F(t,{\bf u},{\bf v})=f_0({\bf u},{\bf v})+f_1({\bf u},{\bf v})\,
(t-\langle {\bf u},{\bf v}\rangle)+R(t,{\bf u},{\bf v})\,H_2(t,{\bf u},{\bf v}),
$$
where $f_0\succeq 0,\; f_1\succeq0$, and $R(t,{\bf u},{\bf v})$ is any continuous function on
$D_m$.




Let $Q=\{q_1,\ldots,q_{n-1}\}\subset{\Bbb S}\sp {n-1}$ with $\rk(Q)=n-1$.  Let $L$ be a linear transformation of the basis $q_1,\ldots,q_{n-1}$ to  an orthonormal basis  $e_1, \ldots, e_{n-1}$  of ${\Bbb R}\sp {n-1}$. Then  $F\in\psd({\Bbb S}\sp {n-1},Q)$ if and only if
$$
F(t,{\bf u},{\bf v})=f_0({\bf u},{\bf v})+f_1({\bf u},{\bf v})\,
(t-\langle L({\bf u}),L({\bf v})\rangle)+R(t,{\bf u},{\bf v})\,H_2(t,L({\bf u}),L({\bf v})),
$$
with $f_0\succeq 0,\; f_1\succeq0$, and any $R(t,L({\bf u}),L({\bf v}))\in C(D_m)$.

In the case $\rk(Q=\{q_1,\ldots,q_m\})\ge n-1$ consider all $Q_i\subset Q$ with $\rk(Q_i)=n-1.$
Let a linear transformation $L_i: {\Bbb R}\sp m\to{\Bbb R}\sp {n-1}$ is defined by $L_i(q)=q$ if $q\in Q_i,\; L_i(q)=0$ if
$q\in Q$ and $\rk(Q_i\cup\{q\})=n$.
Denote by $N$  the number of distinct $Q_i$. It is not hard to prove that $F\in\psd({\Bbb S}\sp {n-1},Q)$ iff 
$$
F(t,{\bf u},{\bf v})=\sum\limits_{i=1}^N{f_i({\bf u},{\bf v})\,F_{i}(t,L_i({\bf u}),L_i({\bf v}))},
$$
where for all $i$:  $f_i({\bf u},{\bf v})\succeq0,\;
 F_{i}(t,L_i({\bf u}),L_i({\bf v}))\in\psd({\Bbb S}\sp {n-1},Q_i).$

\section{Positive definite functions in ${\Bbb R}\sp n$}
Direct extensions of the Bochner - Schoenberg theorem 
and finding bounds on sphere packings in ${\Bbb R}\sp n$ are not so straightforward
because this space is not compact. Different indirect ways
of deriving bounds on sphere packings  in ${\Bbb R}\sp n$ were suggested in the literature
\cite{Kab,gor00,coh03}.

Here we note that the
multivariate Gegenbauer polynomials defined above enable one to
define p.d. functions in ${\Bbb R}\sp n$ as follows:
  $$
    H_k^{(n,m)}(t,x,y,{\bf u},{\bf v}):=(xy)^{k/2}\,
G_k^{(n,m)}(t',{\bf u}',{\bf v}'),
$$
where $0\le m\le n-2,$ $t,x,y\in {\Bbb R},$ ${\bf u},{\bf v}\in {\Bbb R}\sp m$ for $m>0$
and ${\bf u}={\bf v}=0$ for $m=0$, and
   $$
     t'={t}/\sqrt{xy}, \quad {\bf u}'={{\bf u}}/\sqrt{x},
       \quad {{\bf v}}'={{\bf v}}/\sqrt{y}.
$$
The positive semidefiniteness of the polynomials $G_k^{(n,m)}$ (Theorem
3.1) implies the following result.
\begin{theorem}\label{thm:reals}
 Let $e_1, \ldots, e_n  \in {\Bbb R}\sp n$ be an orthonormal basis, and let $p_1,\ldots, p_N$  be points in ${\Bbb R}\sp n$. Then for any $k\ge0$ and  $0\le m \le n-2$  the matrix $\bigl(g_{ij}\bigr)$, where
$$g_{ij}=H_k^{(n,m)}(\langle p_i,p_j\rangle,|p_i|^2,|p_j|^2,p_i,p_j),$$
is positive semidefinite.
\end{theorem}

\begin{remark} It is not hard to find Euclidean analogs of Theorems 3.2, 3.3 and Corollary 3.2
\end{remark}

Theorem 4.1 gives a family of positive-semidefinite constraints for  distance distributions of points in Euclidean spaces.  For instance, consider the simple case of $m=0.$ Now for any matrix $A=\bigl(a_{ij}\bigr)$ of size $N\times N$ we have a matrix
$H_k^{(n)}(A)$ 
which is defined by
$$\Bigl(H^{(n)}_k(A)\Bigr)_{ij}=({a_{ii}a_{jj}})^{k/2}\,G_k^{(n)}\left(a_{ij}/\sqrt{a_{ii}a_{jj}}\right). $$

\begin{cor}
If $A$ is a symmetric positive semidefinite  matrix 
with $\rk(A)\le n$, then for any   positive integer $k$ we have $$H_k^{(n)}(A)\succeq0.$$
\end{cor}
\begin{proof} It is a well known fact: {\em a symmetric matrix $A$ of size $N\times N$ is the Gram matrix of $N$ vectors in ${\Bbb R}\sp n$ if and only if $A\succeq0$ and $\rk(A)\le n$.}
Therefore, there are vectors $p_1,\ldots, p_N$ in ${\Bbb R}\sp n$ such that $A=\bigl(\langle p_i,p_j\rangle\bigr)$. Then
$$\bigl(H_k^{(n)}(A)\bigr)_{ij}=H_k^{(n,0)}(\langle p_i,p_j\rangle,|p_i|^2,|p_j|^2,0,0).$$
Thus Theorem 4.1 yields $H_k^{(n)}(A)\succeq0.$
\end{proof}

For small $k$ it is not hard to give explicit expressions for  $H^{(n)}_k(A)$. Clearly, $H_1^{(n)}(A)=A$.
Since   $G_2^{(n)}(t)=(nt^2-1)/(n-1)$, we have
$$
\bigl(H^{(n)}_2(A)\bigr)_{ij}=\frac{na_{ij}^2-a_{ii}a_{jj}}{n-1}.
$$
Then
$$
H^{(n)}_2(A)=\frac{nA_2-{\bf a}^T{\bf a}}{n-1},
$$
 where
$$
(A_2)_{ij}:=a^2_{ij}, \quad {\bf a}:=(a_{11}, a_{22}, \ldots, a_{NN}).
$$
So 
if $A\succeq0$ and $\rk(A)\le n$, then
$$
nA_2\succeq{\bf a}^T{\bf a}.
$$

\section{Positive semidefinite constraints}
In this section we consider  positive-semidefinite constraints
that are
given by Theorem 3.1.  It is a natural to ask which constraints  are stronger than the others?
We will show that if $m_1>m_2$, then the constraints for $m=m_1$ imply the constraints for $m=m_2$.

Let $T=\bigl(t_{ij}\bigr)$ be a symmetric matrix of size $r\times r$  with $-1\le t_{ij}\le 1,\;  t_{ii}=1,$
and let $U=\bigl(u_{ij}\bigr)$ be a  matrix of size $r\times (n-1)$.
Let $|{\bf u}_i|\le 1$ for all $i=1,\ldots,r$, where
${\bf u}_i:=(u_{i1},\ldots,u_{i,n-1})$. Then we say that a pair $(T,U)_r^n$ is feasible.

\begin{prop} Let $0\le m\le\ell \le n-2, \; d>0$. Let $(T,U)_r^n$ be  a feasible pair.  Suppose
$$\Bigl(G_k^{(n,\ell)}(t_{ij},{\bf u}_i,{\bf u}_j)\Bigr)\succeq 0\; \mbox{ for } \;  k=1,\ldots,d.$$
Then
$$G_k^{n,m}(T,U):=\Bigl(G_k^{(n,m)}(t_{ij},{\bf u}_i,{\bf u}_j)\Bigr)\succeq 0  \; \mbox{ for all } \; k=1,\ldots,d.$$
\end{prop}
\begin{proof} Corollary 2.2 (see also the proof of Corollary 3.1) yields:
$G_k^{n,m}(T,U)$  is a sum of positive semidefinite matrices.
\end{proof}

In Remark 3.4 we considered the case $m=n-1.$ Now we show that the constraints $(3.3), (3.4)$ are strong enough.
\begin{prop} Let $(T,U)_r^n$ be  a feasible pair. Suppose
$$T\succeq\big(\langle {\bf u}_i,{\bf u}_j\rangle\bigr),$$
and for all $i=1,\ldots,r,\; j=1,\ldots,r,$ we have
$$(t_{ij}-\langle {\bf u}_i,{\bf u}_j\rangle)^2=(1-|{\bf u}_i|^2)(1-|{\bf u}_j|^2).
$$
Then there are points $p_1,\ldots,p_r,$ $e_1,\ldots, e_{n-1}$ in ${\Bbb S}\sp {n-1}$ with
$\langle e_i,e_j\rangle=\delta_{ij}$ such that $t_{ij}=\langle p_i,p_j\rangle, \; u_{ik}=\langle p_i,e_k\rangle$ for all $i, j,k$.
\end{prop}
\begin{proof}
 Consider a symmetric matrix $X_0=\bigl(x_{ij}\bigr)$ of size $\ell\times\ell, \; \ell=r+n-1,$ that is defined by:
$x_{ij}=t_{ij} \; \mbox{ for } \; 1\le i\le r, \; 1\le j\le r;$
$
x_{ij}=\delta_{ij}  \; \mbox{ for } \; r< i\le \ell, \; r< j\le \ell;
$
$ x_{ij}=u_{is}, \; s=j-r+n-1,    \; \mbox{ for } \;  1 \le i\le r,\; r< j\le \ell;
$
and $ x_{ij}=x_{ji}    \; \mbox{ for } \;  r< i\le \ell, \; 1\le j\le r.
$

It is easy to see that $(3.4)$ yields $\rk(X_0)\le n.$ On the other hand, it follows from $(3.3)$ that $X_0\succeq0.$ Therefore, there are vectors $q_1,\ldots,q_{r+n-1}$ in ${\Bbb R}\sp n$ such that $X_0=\bigl(\langle q_i,q_j\rangle\bigr).$ Since $x_{ii}=1$, we have $|q_i|=1,$
i.e. $q_i\in {\Bbb S}\sp {n-1}$.
Denote $p_i=q_i$ for $i=1,\ldots,r$, and $e_i=q_{i+r}$ for $i=1,\ldots,n-1.$
\end{proof}

In particular, we obtain that $(3.3), (3.4)$ imply $G_k^{n,m}(T,U)\succeq0$ for all $m: 0\le m\le n-2.$

Actually, $G_k^{n,m}(T,U)\succeq0$ gives constraints only for $T$ and ${\bf u}_i^{(m)}$. For instance, if $m=0$, then we just have $|{\bf u}_i|\le1$. Now using Theorem 3.1 we improve   constraints for $U$.

 Let $e_1, \ldots, e_n$ be an orthonormal basis of ${\Bbb R}\sp n$, and let $p_1,\ldots, p_r$  be points in ${\Bbb S}\sp {n-1}$. Denote $t_{ij}=\langle p_i,p_j\rangle,\; 1\le i\le r,\; 1\le j\le r, \; u_{ik}=\langle p_i,e_k\rangle$, where $1\le k\le n-1.$ Then we have the matrices $T=\bigl(t_{ij}\bigr)$ and $U=\bigl(u_{ij}\bigr)$.

Consider points $\{q_1,\ldots,q_\ell\}=\{p_1,\ldots,p_r,e_{m+1},\ldots,e_{n-1}\}$ in ${\Bbb S}\sp {n-1}$, where $0\le m\le n-2,\; \ell=r+n-m-1.$ Then $x_{ij}=\langle q_i,q_j\rangle$ and $v_{ik}=\langle q_i,e_k\rangle, \; 1\le k\le m,$ define  matrices
$$X_m=X_m(T,U):=\bigl(x_{ij}\bigr),\quad V_m=V_m(T,U):=\bigl(v_{ij}\bigr).$$
Here $V_m$ is well defined for $m>0$. Put $V_0=0.$

Let ${\bf v}_i:=(v_{i1},\ldots,v_{im}),$ where $1\le i\le \ell.\; $
If we apply Theorem 3.1 for the points $\{q_1,\ldots,q_\ell\}$, then for any $k\ge0$ we get
$$
G_k^{n,m}(X_m,V_m):=\Bigl(G_k^{(n,m)}(x_{ij},{\bf v}_i,{\bf v}_j)\Bigr)\succeq 0.
$$

Let $(T,U)_r^n$ be  a feasible pair. Then $X_m(T,U)$ and $V_m(T,U)$ are well defined. It is clear that 
$$
G_k^{n,m}(X_m(T,U),V_m(T,U))\succeq0\; \mbox{ yields } \; G_k^{n,m}(T,U)\succeq0.
$$

Denote by  $\; \Lambda_{d,r}^{n,m},$  $\; 0\le m\le n-2, \; d>0,\; r>0, \; $  the set of all feasible  pairs $(T,U)_r^n$ such that
$$
G_k^{n,m}(X_m(T,U),V_m(T,U))\succeq0 \; \mbox{ for } \; k=1,\ldots,d.
$$

\begin{prop}
$
 \Lambda_{d,r}^{n,n-2}\subset\ldots\subset\Lambda_{d,r}^{n,1}\subset\Lambda_{d,r}^{n,0}.
$
\end{prop}
\begin{proof} It is well known fact: {\em a symmetric matrix $A=\bigl(a_{ij}\bigr)$ of size $\ell\times\ell$ with $a_{\ell\ell}>0$ is positive semidefinite if and only if
$B=\bigl(b_{ij}\bigr)\succeq 0$, where
$$
b_{ij}=a_{ij}-\frac{a_{i\ell}\,a_{j\ell}}{a_{\ell\ell}},\quad 1\le i,j\le\ell-1. \eqno (5.1)
$$}

Let $(T,U)\in \Lambda_{d,r}^{n,m+1}$, where $0\le m\le n-3.$
Denote $$A=G_k^{n,m}(X_m(T,U),V_m(T,U)).$$
Then from the addition theorem (Theorem 2.2) we have:
$$
a_{ij}=\sum\limits_{s=0}^k {h_{s,i,j}}, \quad h_{s,i,j}={C_{k-s}^{n,m+1}({\bf v}_i)\,C_{k-s}^{n,m+1}({\bf v}_j)\,G_s^{(n,m+1)}(x_{ij},{\bf v}_i,{\bf v}_j)}.
$$

Note that $c_{nk0}=1$, where $c_{nks}$ are coefficients in the addition theorem for Gegenbauer polynomials (see subsection 2.1). Indeed, put $\theta_1=\theta_2=\varphi=0$. Then we have $1=G_k^{(n)}(1)=c_{nk0}\,(G_k^{(n)}(1))^2\,G_0^{(n-1)}(1)=c_{nk0}.$ This yields:
$$
h_{0,i,j}=\frac{a_{i\ell}\,a_{j\ell}}{a_{\ell\ell}},\quad 1\le i,j\le\ell-1.
$$

Let $B$ is defined by $(5.1)$. Then $b_{ij}=\sum_{s=1}^k h_{s,i,j}.$ Therefore, $B=\sum_{s=1}^k H_s$, where $H_s=\bigl(h_{s,i,j}\bigr)$. Since $H_s\succeq0$, we have  $B\succeq 0, \; A\succeq0.$
Thus we have proved:
$
\Lambda_{d,r}^{n,m+1}\subset\Lambda_{d,r}^{n,m}.
$
\end{proof}

We denote by $\Delta_{r}^n$ the set of all feasible pairs $(T,U)^n_r$  that satisfy the assumptions of Proposition 5.2 (i.e. $\Delta_{r}^n$ consists of pairs $(T,U)_{r}^n$ that satisfy  $(3.3),(3.4)$).

Note that for $1\le m\le n-2$ we have ${n-1\choose m}$ choices for ${\bf u}^{(m)}$. In other words,
not only $e_1,\ldots,e_m$, but any $m$ vectors from $\{e_1,\ldots,e_{n-1}\}$ can be chosen as a basis of ${\Bbb R}\sp m.$ Denote  by $S\Lambda_{d,r}^{n,m}$ the intersection of all corresponding spaces $\Lambda_{d,r}^{n,m}$.
It is not hard to see that Propositions 5.2, 5.3 imply
\begin{prop}
$
\Delta_{r}^n \subset S\Lambda_{d,r}^{n,n-2}\subset\ldots\subset S\Lambda_{d,r}^{n,1}\subset \Lambda_{d,r}^{n,0}.
$
\end{prop}

\section{Upper bounds for spherical codes}
In this section we set up upper  bounds for spherical codes which are based on multivariate p.d. functions. These bounds extend the famous Delsarte's bound. Note that for the case  $m=1$ this bound is the Bachoc - Vallentin bound \cite{BV3}.

\begin{defn} Consider a vector $J=(j_1,\ldots,j_d)$. Split the set of numbers $\{j_1,\ldots,j_d\}$ into maximal subsets $I_1,\ldots,I_k$ with equal elements. That means,  if $I_r=\{j_{r_1},\ldots,j_{r_s}\}$, then  $j_{r_1}=\ldots=j_{r_s}=a_r$ and all other $j_\ell\ne a_r$.
Without loss of generality it can be assumed that $i_1=|I_1|\ge\ldots\ge i_k=|I_k|>0$.   (Note that we have $i_1+\ldots+i_k=d$.) Denote by $\psi(J)$ the vector $\omega=(i_1,\ldots,i_k)$.

Let
$$
W_d:=\{\omega=(i_1,\ldots,i_k): i_1+\ldots+i_k=d, \; i_1\ge\ldots\ge i_k>0, \; i_1,\ldots,i_k\in {\Bbb Z}\}.
$$
Let $\omega\in W_d$. Denote
$$
\tilde q_\omega(N):=\#\{J=(j_1,\ldots,j_d)\in\{1,\ldots,N\}^d: \psi(J)=\omega\},
$$
$$
q_\omega(N):=\frac{\tilde q_\omega(N)}{N}.
$$
\end{defn}

It is not hard to see that $q_\omega(N)$ is a polynomial of degree $d-1$ for $\omega\in W_d$ and
$$
\sum\limits_{\omega\in W_d} {q_\omega(N)} = N^{d-1}.
$$

\begin{defn}
 For any vector ${\bf x}=\{x_{ij}\}$ with $1\le i<j\le d$ denote by
$A({\bf x})$  a symmetric $d\times d$ matrix $\bigl(a_{ij}\bigr)$  with all $a_{ii}=1$ and  $a_{ji}=a_{ij}=x_{ij}, \; i<j$.

Let   $\; 0<\theta<\pi$ and
$$
X(\theta):=\{{\bf x}=\{x_{ij}\}: x_{ij}\in [-1,\cos{\theta}] \mbox{ or } x_{ij}=1, \; 1\le i<j\le d\}.
$$

Now for any ${\bf x}=\{x_{ij}\}\in X(\theta)$ we  define a vector $J({\bf x})=(j_1,\ldots,j_d)$ such that $j_k=k$ if  there are no $i<k$ with  $x_{ik}=1$, otherwise $j_k=i$, where $i$ is the minimum index with $x_{ik}=1$.

Let $\omega\in W_d$. Denote
$$
D_\omega(\theta):=\{{\bf x}\in X(\theta): \psi(J({\bf x}))=\omega  \; \mbox{ and } \; A({\bf x})\succeq0\}.
$$
Let $f({\bf x})$ be a real function in ${\bf x}$, and let
$$
{B_\omega(\theta,f)}:=\sup\limits_{{\bf x}\in D_\omega(\theta)} {f({\bf x})}.
$$
\end{defn}

Note that the assumption  $A({\bf x})\succeq0$ implies  existence of unit vectors $p_1,\ldots,p_d$ such that $A({\bf x})$ is the Gram matrix of these vectors, i.e. $x_{ij}=\langle p_i,p_j\rangle$. Moreover, if $x_{ij}=1$, then $p_i=p_j$. In particular, $D_{(d)}(\theta)=\{(1,\ldots,1)\}$ and therefore ${B_{(d)}(\theta,f)}=f(1,\ldots,1)$.

\begin{defn}
Let ${\bf x}=\{x_{ij}\}$, where $1\le i<j\le m+2\le n$, and let $A({\bf x})\succeq0$.  
Then there exist $P=\{p_1,\ldots,p_{m+2}\}\subset{\Bbb S}\sp {n-1}$  such that $x_{ij}=\langle p_i,p_j\rangle$.
Let $F({\bf x})$ be a continuous function in ${\bf x}$ with $F({\bf \tilde x}_{k\ell})=F({\bf x})$ for all ${\bf \tilde x}_{k\ell}$ that 
can be obtained by interchanging two points $p_k$ and $p_\ell$ in $P$.  We say that $F({\bf x})\in\psd^n_m$ if for all ${\bf x}$ with  $A({\bf x})\succeq0$ we have  $\tilde F(x_{12},{\bf u_1},{\bf u_2})\in\psd({\Bbb S}\sp {n-1},Q({\bf x}))$, where ${\bf u_i}=(x_{i3},\ldots,x_{i,m+2})$, $Q({\bf x})=\{p_3,\ldots,p_{m+2}\}$, and $\tilde F(x_{12},{\bf u_1},{\bf u_2})=F({\bf x})$.
\end{defn}


For the classical case $m=0$  Schoenberg's theorem says that $f\in \psd^n_0$ if and only if
$f(t)=\sum_{k}{f_kG_k^{(n)}(t)}$ with all $f_k\ge 0$.
Theorem 3.3 (see also \cite{BV,BV2,BV3}) yields $F({\bf x})\in\psd^n_1$ if and only if
$$
F(x_{12},x_{13},x_{23})=\sum\limits_{k}{f_k(x_{12},x_{13},x_{23})\,G_{k}^{(n,1)}(x_{12},x_{13},x_{23})},
$$
where $f_k\succeq0$  for all $k,$ and
$$F(x_{12},x_{13},x_{23})=F(x_{13},x_{12},x_{23})=F(x_{23},x_{13},x_{12}).$$
Using Corollary 3.2 it is possible to describe the class of functions $\psd^n_m$, $m\ge 2$.

\medskip

Let $C$ be an $N$-element subset of the unit sphere ${\Bbb S}\sp {n-1}\subset {\Bbb R}\sp n$.  It is called an    $(n,N,\theta)$ {\it spherical code} if every pair  of distinct points  $(c,c')$ of $C$ have inner product $\langle c,c'\rangle$ at most $\cos{\theta}$.


\begin{theorem}
Let $f_0>0, \; 0\le m\le n-2$, and $F({\bf x})=f({\bf x})-f_0\in\psd^n_m$. 
Then an $(n,N,\theta)$ spherical code satisfies
$$
f_0N^{m+1}\le \sum\limits_{\omega\in W_{m+2}} {B_\omega(\theta,f)\,q_\omega(N)}.
$$
\end{theorem}

\begin{proof} Let $C$ be an $(n,N,\theta)$ spherical code. Define
$$
S=\sum\limits_{{\bf c}=(c_1,\ldots,c_d)\in C^d} {f(\{\langle c_i,c_j\rangle\})}, \quad d=m+2.
$$
Then
$$
S= \sum\limits_{\omega\in W_{d}}\;\sum\limits_{{\bf c}: \psi({\bf c})=\omega}  {f(\{\langle c_i,c_j\rangle\})}\le \sum\limits_{\omega\in W_{d}} {B_\omega(\theta,f)\,\tilde q_\omega(N)}.
$$
On the other hand, since  $F\in\psd^n_m$ we have
$$
\sum\limits_{{\bf c}\in C^d} {F(\{\langle c_i,c_j\rangle\})}\ge 0.
$$
Thus
$$
S=\sum\limits_{{\bf c}\in C^d} {(f_0+F(\{\langle c_i,c_j\rangle\}))}\ge f_0\,N^d.
$$
\end{proof}

It is easy to see for $m=0$ that $q_{(2)}(N)=1, \; q_{(1,1)}(N)=N-1,$ and  $B_{(2)}(\theta,f)=f(1)$. Therefore, from Theorem 6.1  we have
$$f_0N\le f(1)+B_{(1,1)}(\theta,f)(N-1).$$ Suppose $B_{(1,1)}(\theta,f)\le0$, i.e. $f(t)\le 0$ for all $t\in [-1,\cos{\theta}]$. Thus for $(n,N,\theta)$ spherical code we obtain
$$
N\le \frac{f(1)}{f_0}.
$$
This upper bound is called Delsarte's bound.

\medskip

\medskip

The Bachoc-Vallentin bound \cite[Theorem 4.1]{BV3}
is the bound in Theorem 6.1 for $m=1$ and $B_{(1,1,1)}(\theta,f)\le0$. Indeed, let $B_{(2,1)}(\theta,f)\le B$. Since $q_{(3)}(N)=1, \; q_{(2,1)}(N)=3(N-1),$  and $B_{(3)}(\theta,f)=f(1,1,1)$, we have
$$
f_0N^2\le f(1,1,1)+3(N-1)B.
$$

Let us consider Theorem 6.1 also for the case $m=2$ with $B_{(1,1,1,1)}(\theta,f)\le0$. Let $B_{(3,1)}(\theta,f)\le B_1, \; B_{(2,2)}(\theta,f)\le B_2,$ and $B_{(2,1,1)}(\theta,f)\le B_3$. Then
$$
f_0N^3\le f(1,1,1,1,1,1)+4(N-1)B_1+3(N-1)B_2+6(N-1)(N-2)B_3.
$$

Let $f({\bf x})$ be a polynomial of degree  $d$.
Then the assumptions in Theorem 6.1
 can be written as  positive semidefinite constraints for the coefficients of $F$ (see for details \cite{BV,BV2,BV3,GST,Schr}).
Actually, the bound given by Theorem 6.1 can be obtained as a solution  of an SDP (semidefinite programming) optimization problem. In \cite {BV,BV2} using numerical solutions of the SDP problem for the case $m=1$
has obtained new upper bounds for the kissing numbers and for the one-sided kissing numbers in several dimensions $n\le 10$.

However, the dimension of the corresponding SDP problem is growth so fast
whenever $d$ and $m$ are increasing
that this problem can be treated numerically only for relatively small $d$ and small $m$.
It is an interesting problem to find  (explicitly) 
suitable  polynomials $F$ for Theorem 6.1 and using it to obtain new bounds for spherical codes.





\end{document}